\documentstyle[12pt]{article}
\begin{document}
\begin{center}
\Large
STRUCTURE  THEORY  FOR  ONE  CLASS  OF  LOCALLY
FINITE  LIE  ALGEBRAS.
\vspace{.8cm}

\large
L.A. Simonian
\footnote{Send correspondence to the author at \emph{levsimonian@hotmail.com} or 1060 Ocean Avenue F8, Brooklyn, NY 11226, USA}

\end{center}

\begin{center}

\large{Abstract}

\end{center}
In this paper I consider locally finite Lie algebras of
characteristic zero satisfying the condition that for every finite
number of elements $x_{1}, x_{2},\cdots ,x_{k}$ of such an algebra
$L$ there is finite-dimensional subalgebra $A$ which contains
these elements and $L(adA)^n\subset A$ for some integer $n$. For
such algebras I prove several structure theorems that can be
regarded as generalizations of the classical structure theorems of
the finite-dimensional Lie algebras theory.
\begin{center}
---------------------
\end{center}

\paragraph{INTRODUCTION.}

The subject of this article is similar to that of Refs.\cite{1,2,3}. I consider
locally finite Lie algebras of characteristic zero. A Lie algebra
is called locally finite if every its finite subset is contained
in a finite-dimensional subalgebra. We will study representations
$R = ( M , L )$ of such algebras. A representation $R = ( M , L )$
is a homomorphism of a Lie algebra $L$  into the algebra of
linear transformations of a linear space  $M$.  We will assume that $R = ( M , L  )$ satisfies the following condition:

{\sl (1) for every  $x\in L$  the linear transformation  $x^R$  has
the Fitting null component of a finite co-dimension}.
Here  $x^R$ is the linear
transformation that corresponds  to  $x\in L$  in the representation $R$.

For such linear transformation the following is true:

{\sc LEMMA 1}. {\sl Let  $A$  be a linear transformation of a
linear space  $M$ of an infinite dimension,  $M$ be the Fitting
null component of  $A$, that is, the set of all  $x\in  M$  with
$xA^m = 0$  for some integer  $m$. If $M /M_0$  has finite
dimension,  then  $M$  is a direct sum  $M_0\oplus
M_1$, where $M_1$ is a finite-dimensional space invariant under $A$, and
the transformation induced by $A$  in $M_1$  is an automorphism}.

Lemma 1 allows us to determine the trace  $trA$  of $A$  as a
trace of the linear transformation, which is induced by  $A$  in
the finite-dimensional space  $M /M_0$. Besides, given a
representation  $R = ( M , L )$  satisfying condition  (1),  this
Lemma allows us to construct a decomposition
$$M  =  M_{\rho}\oplus M_{\sigma}\oplus\ldots\oplus M_{\tau}
\oplus M_0$$ into weight spaces relative to a nilpotent subalgebra
$H$  of  $L$. This decomposition has the same properties as  in
the finite-dimensional case and is used in the proof of

{\sc THEOREM 3}. {\sl Let  $L$  be a locally finite Lie algebra of
characteristic zero. Suppose  $L$  has a representation  $R = (M , L)$  in a
vector space  $M$  satisfying the condition  (1)  and the associative span
$\cal A$ of  $L$  in the representation  $R = ( M , L )$ is locally finite.
If the kernel of the representation  $R$  is locally solvable and
$tr(x^R)^2 = 0$  for every  $x\in L'$,  then  $L$  is locally solvable}.

Theorem 3 can be considered as a generalization of the Cartan's
criterion for solvability.

Subsequent results are obtained for Lie algebras and their
representations that satisfy the following conditions:

{\sl (2) representations  $R = ( M , L )$: for every
finite-dimensional subagebra $A$ of $L$ there is such an integer
$n$ that $MA^n$ has finite dimension.}

{\sl (3) Lie algebra  $L$: for every finite set
$x_{1},x_{2},\cdots,x_{k}\in L$ there is a finite-dimensional
subagebra $A$ which contains these elements and for which
$L(adA)^n\subset A$ for some integer $n$.}

{\sc COROLLARY}. {\sl An algebra  $L$  satisfying the condition  (3)  is
locally solvable if and  only if  $tr (adx)^2 = 0$  for every  $x\in L'$}.

Since the intersection of a finite number of finite-codimensional
subspaces has finite codimension, the trace can be defined simultaneously for any finite number of transformations for algebras that satisfy condition
(2). All properties of the usual finite-dimensional
trace are true in this case. Therefore, for the representation $R =
( M , L )$ satisfying the condition (2) we can define a trace form
$$f (a,b) = tr a^Rb^R,  a, b \in L.$$
In particular, if $L$ satisfies  the condition  (3),  we obtain the form
$$K (a,b) = tr (ada)(adb)$$ It is natural to name this form the Killing 
form of  $L$.

Lie algebra  $L$  is said to be semi-simple if its locally solvable ideal
equals to 0, that is, if $L$  has no non-zero locally solvable ideals.

{\sc THEOREM 4}. {\sl Let  $L$  be a locally finite semi-simple Lie algebra
of characteristic 0, and   $R = ( M , L )$  be an arbitrary faithful
representation satisfying the condition  (2).  Then the trace form  $f
(a,b)$ is non-degenerate.

If  $L$  satisfies the condition  (3)  and the Killing form
$K (a,b)$ is non-degenerate, then  $L$  is semi-simple.}

This theorem can be regarded as a generalization of the Cartan's
criterion for semi-simplicity. Finally, the following theorem
generalizes the structure theorem.

{\sc THEOREM 6} (Structure Theorem). {\sl Let  $L$  be a
semi-simple Lie algebra that satisfies condition (3). Then
$L$ is a subdirect sum of a set of finite-dimensional simple
algebras.}

\paragraph{1. WEIGHT  SPACES.}

{\sc PROOF OF  LEMMA 1}. Let  $N$  be a finite-dimensional subspace of  $M$
such that $M = N\oplus M_0$  and let  $e_1, e_2,\ldots , e_n$  be a basis of
the $N$. We have  $e_iA=\sum\limits_{j=1}^na_{ij}e_j + m_i$, where
$m_i\in M_0$.  For any  $m_i$  there is an integer number  $n_i$  such that
$m_iA^{n_i} =0$.
Let  $P$  denote the linear subspace spanned by  $m_i$, $m_iA$,\ldots ,
$m_iA^{n_1-1} , i=1, 2,\ldots , n$.  It is clear that  $P$  is invariant
under $A$ and has a finite dimension. From the equality
$e_iA=\sum\limits_{j=1}^na_{ij}e_j +m_i$ it follows that  $NA\subset P + N$.
Hence  $K = P + N$   is invariant under $A$.  We have $M = K + M_0$.  Since
$K$ is finite-dimensional, it can be represented as  $K = K_1\oplus K_0$,
where $K_1$ and  $K_0$  are, respectively, the Fitting one and the Fitting
null components of $K$ relative to the transformation induced by $A$ in $K$.
Then $M = K_1 + K_0 + M_0$. Since $K_0\subset M_0$, then $M = K_1 + M_0$. We
shall show now  that this sum is direct. Since $K$ is finite-dimensional
there is $t$ such that $yA^t=0$ for any $y\in K\cap M_0$. On the other hand,
there exists $s$ such that $K_1=KA^s=KA^{s+1}=\cdots$. Let $r=max(s, t)$.
Then $K_1=KA^r$ and $yA^r=0$ for any $y\in K\cap M_0$. Now let $x\in K_
1\cap M_0$. Then from the equality $K_1=KA^r$  it follows that $x=yA^r$ for
some $y\in K$. On the other hand, since $x\in K_ 1\cap M_0$, it holds
$xA^r=0$. But then $0 = xA^r = yA^{2r}$.  Hence $y \in K\cap M_0$ and,
consequently, $yA^r=0$.  But then $x=yA^r$ is equal to zero and $K_1\cap
M_0=0$. From the construction of $K_1$ it follows, that $A$ is a linear
transformation acting in $K_1$ as an automorphism. Therefore, the only possibility is to put $M_1=K_1$.

This lemma may be considered as a generalization of the well-known Fitting's
lemma. Let $M_0$ and $M_1$ are the Fitting null and
the Fitting one components of $M$ relative to $A$, respectively.

As in [4], a linear transformation $A$ will be called algebraic, if every
vector $x\in M$ is contained in a finite-dimensional subspace that is
invariant under~$A$.

{\sc LEMMA 2}. {\sl If the Fitting null component $M_0$ of $M$ relative to
$A$ has finite codimension, then $A$ is an algebraic linear transformation.}

{\sc PROOF}. By Lemma 1, $M=M_1\oplus M_0$, where $M_1$ has finite dimension
and A acts in $M_1$ as an automorphism. Let $x\in M$. We shall show that
the dimension of the smallest subspace that contains $x$ and is invariant relative to $A$ is finite. Let $x=y + z$, where $y \in M_1$ and $z\in M_0$.
Since $M_1$ is invariant under $A$, $yA\in M_1$. On the other hand, $zA^m= 0$ for some integer $m$.  Let $N$ be the subspace, which is
generated by  $M_1$, $z,zA,\ldots,zA^{m - 1}$. It is clear that $x\in N$,
$N$ is invariant under $A$, and $N$ has finite dimension. The lemma is proved.

Let the characteristic roots of $A$ be in the base field $\Phi$ and let
$M_1=M_\alpha\oplus M_\beta\oplus\cdots\oplus M_\gamma$
be the decomposition of $M_1$ into the weight spaces relative to $A$. Then
$M=M_0\oplus(\bigoplus\limits_{\alpha\ne0}M_\alpha)$.  Also we have that all
$M_\alpha$  with $\alpha\ne0$ are of finite dimension. It is worth
recalling that by definition $x\in M_\alpha$, if and only if $x(A -\alpha
E)^m=0$  for some integer m.

We will require some known results which I outline here for completeness. Let $\cal A$ be an associative algebra, and $a\in\cal A$. Let
us consider the inner derivation $D_a:x\to x'= [x, a]$ in $\cal A$. If we
denote $x^{(k)}=(x^{(k - 1)})', x^{(0)}=x$, then the following formulas
hold:
$$xa^k = a^kx+{k\choose1}a^{k-1}x'+{k\choose2}a^{k -2}x''+\cdots+x^{(k)}$$
$$a^kx = a^kx-{k\choose1}x'a^{k-1}+{k\choose2}x''a^{k-2}+\cdots+(- 1)^k
x^{(k)}$$
$$x\phi(a)=\phi(a)x+\phi_1(a)x'+\phi_2(a)x''+\cdots+x^{(r)},$$
where $\phi(\lambda)$ is a polynomial of degree $r$ and
$\phi_k(\lambda) =\phi^{(k)}(\lambda)/k!$.

{\sc LEMMA 3} \cite{6}. {\sl Let $A,B$  be linear transformations in a
vector space $M$ satisfying $B(adA)^u=0$ for some integer $u$. Let
$\mu(\lambda)$  be a polynomial and let $M_{\mu
A}=\{x|\mbox {$x\mu(A)^m=0$ for some integer $m$}\}$. Then
$M_{\mu A}$ is invariant under $B$}.

{\sc PROOF}. Let $x\in M_{\mu A}$ and suppose that $x\mu(A)^m=0$.
Putting $\phi(\lambda)=\mu(\lambda)^{mu}$, we obtain $B\phi(A)=\phi(A)B+
\phi_1(A)B'+\cdots+\phi_{u-1}(A)B^{(u - 1)}$.
Since $\phi_0(\lambda)=\phi(\lambda),\phi_1(\lambda),\ldots,\phi_{u -1}(
\lambda)$  are divisible by $\mu(\lambda)^m$, $x\phi_j(A)=0$, $0\le j \le u
- 1$. Therefore $xB\phi(A)=0$ and $xB\in M_{\mu A}$.

{\sc COROLLARY}. {\sl If $B(adA)^u=0$ then the weight spaces $M_{\alpha}$
are invariant under B}.{\sloppy

}Let $R=(M , L)$ be a representation of a Lie algebra $L$ in a vector
space $M$ of infinite dimension, satisfying the condition (1) and let
characteristic roots of every $A^R$, $A\in L$,  lie in the base field.

We shall also assume that for every finite-dimensional subalgebra $H$ of $L$
there is such  an integer $m$ that $L(adH)^m=0$.

{\sc THEOREM 1}. {\sl If $H$ is an finite-dimensional subalgebra of $L$,
then $M$ can be decomposed as
$(\bigoplus\limits_{\alpha\ne0}M_{\alpha})\oplus M_{0H}$ where $M_{\alpha}$,
$\alpha\ne 0$, are finite-dimensional weight spaces relative to $L$ with the
weights $\alpha$, and $M_{0H}$ is a weight space relative to $H$ with the
weight $\alpha=0$. The dimension of $M_{0H}$ is infinite}.

{\sc PROOF}. First let's show that for every $x\in M$ the smallest subspace
$N$ that is invariant under $H$ and contains $x$ has finite dimension. Indeed
(see also \cite{4}), if $B_1,B_2,\ldots,B_r$ is a basis of $H$, then
$N$ is the linear span of the set of all elements of the form
$xB_{i_1}^{m_1}B_{i_2}^{m_2}\ldots B_{i_k}^{m_k}$, $i_1\le
i_2\le\ldots\le i_k$. Since all $B_i$ are algebraic, there is only a finite
set of linearly independent elements of a given form and all of them may be
found among the elements $xB_{i_1}^{m_1}B_{i_2}^{m_2}\ldots B_{i_k}^{m_k}$
for which $m_j\le s_j$, where $s_j$  are  integers and $j=1, 2,\ldots,k$.
Suppose now that every element $A$ of a finite subset $F\subset H$
is locally nilpotent, that is, for $A$ the following condition is
satisfied: for any $x\in M$ there exists $m$ such that $xA^m=0$.  We shall
show that this condition holds for every element of the subalgebra $\{F\}$
generated by the set $F$.  Indeed, since $H$ is nilpotent, $F$ is contained
in Jacobson radical of representation $({N}, H)$
\cite{5}.  Hence  $\{F\}$ is contained in it. But this means that $\{F\}$
consists of nilpotent relative to $N$ transformations and hence the
given condition holds for the elements of  $\{F\}$.

Since $H$ is nilpotent, there exists a chain of ideals
$0\subset H_1\subset H_2\subset\ldots H_{n-1}\subset H_n = H$  such that
$\dim H_{i+1}/H_i=1$. Take an arbitrary element $A\in H_1$ and let $M =
\bigoplus\limits_{\alpha_A}M_{\alpha_A}$ be a decomposition of $M$ into the
weight spaces relative to $A^R$. From Lemma 3  it follows that all
$M_{\alpha_A}$ are invariant under $L$. Besides there is just a finite
number of the subspaces $M_{\alpha_A}$, and $\dim M_{\alpha_A}<\infty$ if
$\alpha_A\ne0$. Therefore
$M_{1A}=\bigoplus\limits_{\alpha_A\ne0}M_{\alpha_A}$
can be decomposed into a direct sum of a finite number of weight spaces
relative to $L$ (see \cite[p. 43]{6}). We recall that a map
$\alpha:A\to\alpha(A)$ of $L$ into the base field $\Phi$ is called the
weight of $M$ relative to $L$ if there exists a nonzero element
$x\in M$ such that $x(A^R - a(A))^m = 0$ for all $A\in L$. Here  $m$ is an
integer which depends on $x$ and $A$. The set of elements (zero included)
satisfying this condition forms the subspace that is called the weight
subspace. It should be recalled that $M_{0A}$ also is invariant under $L$.
Let $B$ be an element of  $H_2\setminus H_1$ and
$M_{0A}=\bigoplus\limits_{\alpha_B}M_{\alpha_B}$ be a decomposition of
$M_{0A}$ into the weight  spaces relative to $B^R$. The subspace $M'_{1B}=
\bigoplus\limits_{\alpha_B\ne0}M_{\alpha_B}$ has a finite dimension, is
invariant under $L$ and can be decomposed into a direct sum of a finite
number of weight spaces relative to $L$. If we add this decomposition to
decomposition of $M_{1A}$, we obtain that $M=(\bigoplus\limits_{\alpha\ne0}M
_\alpha)\oplus M'_{0B}$. $A^R$ and  $B^R$ act in $M'_{0B}$ as locally
nilpotent transformations. Therefore, the subalgebra $\{A, B\}$ consists of
the locally nilpotent in $M'_{0B}$ transformations \cite{7}. Continuing in
this way we obtain - by virtue of finite dimensionality of $H$ - the
statement of the theorem.

If $L$ is nilpotent we may combine Lie's theorem  with Theorem 1 to obtain the following

{\sc THEOREM 2}. {\sl If $L$ is finite-dimensional, then $M$ is a direct sum
of weight spaces $M_{\alpha}$, and the matrices in the weight space
$M_{\alpha}, \alpha\ne0$ can be taken simultaneously in the form
$$A_\alpha=\left(\begin{array}{cccc}\alpha(A) & 0 & \ldots & 0
\\ ** & \alpha(A) & \ldots & 0 \\
\ldots & \ldots & \ldots & \ldots \\ **
& * & \ldots & \alpha(A) \end{array}\right)$$}

This theorem is proved in exactly the same way as in Ref.\cite{6}. In a similar
fashion we obtain

{\sc COROLLARY}. {\sl The weights  $\alpha: A\to\alpha(A)$ are linear
functions on $L$ which vanish on $L'$}.

\paragraph{2. CARTAN'S  CRITERION.}

Let $L$ be a finite-dimensional Lie algebra, $H$ be a nilpotent subalgebra
of $L$, $R=( M ,  L )$ be a representation of $L$ in a vector space $M$
satisfying the condition (1).

{\sc PROPOSITION 1}. {\sl Let
$$M=M_\rho\oplus M_\sigma\cdots M_\tau\oplus M_0$$
$$L=L_\alpha\oplus L_\beta\oplus\cdots\oplus L_\gamma\oplus L_0$$
be the decompositions of $M$ and $L$ into weight spaces relative to $H$.
(The existence of the first decomposition was proved in the previous
section). Then  $M_\rho L_\alpha\subset M_{\rho+\alpha}$ if  $\rho+\alpha$
is the weight of $M$ relative to $H$; otherwise $M_\rho L_\alpha=0$}.

{\sc PROOF}. For every $x\in M$ and $A,B\in L$ we have the equality
$xA(B-\rho I-\alpha I)=x(B-\rho I)A+x(A(adB-\alpha I))$. If $x(B-\rho I)^m
= 0$ and $A(adB -\alpha I)^n=0$ then by repeating this equality we obtain $xA
(B-\rho I-\alpha I)^{m+n+1}=0$. Here
$\rho=\rho(B)$, $\alpha=\alpha(B)$, and $I$ is an identity operator of $M$.

It is also true that $[L_\alpha, L_\beta]\subset L_{\alpha+\beta}$ if
$\alpha+\beta$ is a root of $L$ and $[L_\alpha,L_\beta]=0$ otherwise
(see \cite[p. 64]{6}).

Suppose now that $H$ is a Cartan subalgebra. Then $H=L_0$,  the root
module corresponding to the root  0.  Also, we have $L'=[L,
L]=\sum\limits_{\alpha,\beta}[L_\alpha,L_\beta]$, where the sum is taken
over all roots  $\alpha,\beta$, and $L'\cap
H=\sum\limits_\alpha[L_\alpha,L_{-\alpha}]$,
where the summation is taken over all $\alpha$ such that $-\alpha$ is also
a root (see \cite[p. 67]{6}).

Let $A$ be a linear transformation of $M$ with the Fitting null component
$M_{0A}$ of a finite codimension. Then, as noted in Introduction, $trA$
can be defined as the trace of the linear transformation induced by $A$ in
the quotient space $M/M_{0A}$. It is easy to see that $trA$ is equal to
the trace of the linear transformation induced by $A$ in the quotient space
of $M$ by any invariant under $A$ subspace of a finite codimension
contained in $M_{0A}$.

{\sc LEMMA 4}. {\sl Let $\Phi$ be algebraically closed of characteristic
0. Under the assumptions of this section let $H$ be a Cartan subalgebra of
$L$ and let $\alpha$ be a root such that $-\alpha$ is also a root. Let
$e_\alpha\in L_\alpha$, $e_{-\alpha}\in L_{-\alpha}$,
$h_\alpha=[e_\alpha,e_{-\alpha}]$. Then $r(h_\alpha)$ is a rational multiple
of $\alpha(h_\alpha)$ for every weight $\rho$ of $H$ in $M$}.

{\sc PROOF}. Let $M_0^\alpha=M_0+\sum\limits_iM_{i\alpha}$, $i=0,\pm1,
\pm2,\ldots$. Let's turn to the quotient space $\overline M=M/M_0^\alpha$.
If $M_0^\alpha$ is invariant under $x\in L$, then the operator induced  by
x in $\overline M$ is denoted as $x^{\overline R}$. Consider functions
of the form $\rho(h)+i\alpha(h), i = 0,\pm1,\pm2,\ldots$, which are weights,
and form the subspace $N=\sum\limits_i\overline M_{\rho+i\alpha}$   where
$\overline M_{\rho+i\alpha}=M_{\rho+i\alpha}+M_0^\alpha/M_0^\alpha$  and the
sum is taken over the corresponding weight spaces of the representation
$R=(M,L)$. $N$ is invariant relative to $H$ and, by Proposition 1, it is
also invariant relative to the linear transformations  $e_\alpha^{\overline
R}$ and $e_{-\alpha}^{\overline R}$. Thus, if $tr_N$ denotes the trace of an
induced mapping in $N$, then $tr_Nh_\alpha^{\overline
R}=tr_N[e_\alpha^{\overline R},e_{-\alpha}^{\overline R}]=0$. On the other
hand, the restriction of $h_\alpha^{\overline R}$ to $\overline
M_\sigma=M_\sigma+M_0^\alpha/M_0^\alpha$ has the single characteristic root
$\sigma(h_\alpha)$. Hence $0=tr_Nh_\alpha^{\overline
R}=\sum\limits_in_{\rho +i\alpha}(\rho+i\alpha)(h_\alpha)$ where
$n_{\rho+i\alpha}= \dim M_{\rho+i\alpha}$. Thus we have
$(\sum\limits_in_{\rho+i\alpha})\rho(h_\alpha)+(\sum\limits_iin_
{\rho+i\alpha})\alpha(h_\alpha)=0$.  Since $\sum\limits_in_{\rho+i\alpha}$
is a positive integer, this shows
that $\rho(h_\alpha)$ is a rational multiple of $\alpha(h_\alpha)$.

{\sc PROOF OF  THEOREM 3}. Assume first that the base field $\Phi$ is
algebraically closed. It suffices to prove that $C'\subset C$ for every
finite-dimensional subalgebra $C$ of the algebra $L$.  Hence we shall have
that $C\supset C'\supset C''\supset C^{(k)}=0$. We therefore suppose
that there exists a finite-dimensional subalgebra $C$ such that $C'=C$.
Let $H$ be a Cartan subalgebra of $C$ and let
$$M=M_\rho\oplus M_\sigma\oplus\cdots\oplus M_\tau\oplus M_0$$
$$C=C_\alpha\oplus C_\beta\oplus\cdots\oplus C_\gamma\oplus C_0.$$
be the decomposition of $M$ and $C$ into weight spaces relative to $H$.
Then the formula
$H\cap C'=\sum\limits_\alpha[C_\alpha,C_{-\alpha}]$ implies that
$H=\sum\limits_\alpha[C_\alpha,C_{-\alpha}]$ summed on $\alpha$ such that
$-\alpha$ is also a root. Choose such $\alpha$, let  $e_\alpha\in C_\alpha$,
$e_{-\alpha}\in C_{-\alpha}$, and consider the element $h_\alpha=[e_\alpha,
e_{-\alpha}]$. The formula $H=\sum\limits_\alpha[C_\alpha,C_{-\alpha}]$
implies that every element of $H$ is a sum of terms of the form $[e_\alpha,
e_{-\alpha}]$. Let us turn to the quotient space $\overline M=M/M_0$ and
denote by $h_\alpha^{\overline R}$ the operator induced by $h_\alpha$ in
$\overline M$. The restriction of $h_\alpha^{\overline R}$ to
$\overline M_\rho=M_\rho+M_0/M_0$ has the single characteristic root
$\rho(h_\alpha)$. Hence the restriction of $(h_\alpha^{\overline R})^2$  has
the single characteristic root $\rho(h_\alpha)^2$. Let $n_\rho$ be the
dimension of $M_\rho$. Then we have  $tr(h_\alpha^{\overline
R}))^2=\sum\limits_{\rho\ne0}n_\rho(\rho(h_\alpha))^2$. On the other hand,
$tr(h_\alpha^{\overline R})=tr(h_\alpha^R)$. Thus
$\sum\limits_{\rho\ne0}n_\rho(\rho(h_\alpha))^2=0$ since $tr(h_\alpha^R)^2 =
0$. By the Lemma 4, $\rho(h_\alpha)=r_\rho\alpha(h_\alpha)$, where $r_\rho$
is rational. Hence $\alpha(h_\alpha)^2(\sum\limits_{\rho\ne0}n_\rho
r_\rho^2)=0$. Since $n_\rho$ are positive integers, this implies that
$\alpha(h_\alpha)=0$ and $\rho(h_\alpha)=0$.  Since $\rho$ are linear
functions and every $h\in H$ is a sum of elements of the form $h_\alpha,
h_\beta,\ldots$,  we see that $\rho(h)=0$. Thus $0$ is the only weight for
$M$, that is, we have $M=M_0$. If $\alpha$ is a root, then the condition
$$M_\rho C_\alpha=\left\{\begin{array}{lcc}0\qquad\mbox{if $\rho+\alpha$ is
not a weight of $M$}&\\ \subset M_{\rho+\alpha}\qquad\mbox{if $\rho+\alpha$
is a weight}&\end{array}\right.$$
implies that $MC_\alpha=0$ for every $\alpha\ne0$.  Hence $C_\alpha\oplus
C_\beta\oplus\cdots\oplus C_\gamma$, $\alpha, \beta,\ldots, \gamma\ne0$
is contained in the kernel $K$ of representation $(M,C)$.
Hence  $C/K$  is a homomorphic image of $H$.  It follows that the $C/K$ is
nilpotent. According to our assumptions the kernel  K  is solvable, and it
follows that $C$ is solvable which contradicts $C'=C$.

If the base field is not algebraically closed, then let $\Omega$  be its
algebraic closure. Then  $(M_\Omega, L_\Omega)=R_\Omega$ is the
representation of $L_\Omega$ in $M_\Omega$ and $K_\Omega$ is the kernel of
this representation if $K$ is the kernel of $R=(M,L)$. Since $K$ is
locally solvable, $K_\Omega$ is locally solvable. Next we note that
$tr(x^R)^2=0$ and $trx^Ry^R=try^Rx^R$ imply that
$trx^Ry^R=(1/2)tr(x^Ry^R+y^Rx^R)
=(1/2)(tr(x^R+y^R)^2-tr(x^R)^2-tr(y^R)^2)=0$.
Hence if $x_i\in L$  and $\omega_i\in\Omega$, then
$tr(\sum\limits_i\omega_ix_i^R)^2=\sum\limits_i\omega_i\omega_jtrx_i^Rx_j^R=0$.
To prove that the condition  (1)  holds we use the fact that associative
span $\cal A$ of $L$ in the representation   $R=(M,L)$ is locally finite. We
need to show that $(\sum\limits_{i=1}^m\omega_ix_i)^R$  has the Fitting null
component of finite codimension for every  $x_1,x_2,\ldots,x_m\in L$  and
$\omega_1,\omega_2,\ldots,\omega_m\in\Omega$. Since $\cal A$ is locally
finite, the subalgebra $A$ in $\cal A$ generated by $x_1,x_2,\ldots,x_m$
has a finite dimension. Therefore there exists an integer $u$ such that  $y(
x^R)^u=0$. Here $x$ is an arbitrary linear combination of $x_1,x_2,\ldots,
x_m$ with coefficients from the base field $\Phi$, and $y$ is an arbitrary
element from the Fitting null component $M_{0X}$ of the linear
transformation $x^R$.

Take $(u + 1)^m$  elements of the form
$$k_{i_1}x_1+m_{i_2}x_2+n_{i_3}x_3+\cdots+s_{i_m}x_m$$, where
$i_1,i_2,i_3,\ldots,i_m$  receive their values  $1, 2, 3,\ldots, u+1$
independently, and $k_1,k_2,\ldots,k_{u+1},m_1,m_2,\ldots,m_{u+1},n_1,n_2
,\ldots,n_{u+1},\ldots,s_1,s_2,\ldots,s_{u+1}$ are pairwise different
nonzero integers. The intersection $N$ of the Fitting null components
of these elements has a finite codimension. Let us prove that for any $y\in
N$ the equality $y(\omega_1x_1+\omega_2x_2+\cdots+\omega_mx_m)^u=0$
takes place for every $\omega_1,\omega_2,\ldots,\omega _m\in\Omega$. We have
$$(k_jx_1+m_{i_2}x_2+n_{i_3}x_3+\cdots+s_{i_m}x_m)^u=0, j=1,
2,\ldots,u+1$$ for any $m_{i_2},n_{i_3},\ldots,s_{i_m}$ which are taken from
the set of integers shown above. It follows immediately that
$$yP_0+k_jyP_1+k_j^2yP_2+\cdots+k_j^uyP_u=0,j=1,2,\ldots,u+1$$ where $P_i=
P_i(x_1,m_{i_2}x_2,n_{i_3}x_3,\ldots,s_{i_m}x_m)$ is the homogeneous
component of $i$-th degree relative to $x_1$ of 
$(k_jx_1+m_{i_2}x_2+n_{i3}x_3
+\cdots+s_{i_m}x_m)^u$. Since the determinant of this system is
Vandermonde's determinant, it follows that
$yP_0=yP_1=yP_2=\cdots=yP_u=0$. Next for $P_i,i= 0, 1, 2,\ldots,u$,  we have
$yP_i(x_1,m_lx_2,n_{i_3}x_3,\ldots,s_{i_m}x_m)=0,l=1,2,\ldots,u-i+1$.
It immediatelly follows that
$$yP_{i0}+m_lyP_{i1}+m_l^2yP_{i2}+\ldots+m_l^{u
-i}yP_{i(u-i)}=0,l=1,2,\ldots,u-i+1$$ where
$P_{ij}(x_1,x_2,n_{i3}x_3,\ldots,s_{i_m}x_m)$ is a component of
$P_i$ which is homogeneous of degree $i$ relative to $x_1$ and of degree $j$
relative to $x_2$. Since the determinant of this system is Vandermonde's
determinant, it follows that $yP_{i0}=yP_{i1}=\cdots=yP_{i(u -i)}=0$.

Continuing in this way we obtain that $yP_{i_1i_2\ldots i_m}=0,
i_1,i_2,\ldots,i_m=1,2,\ldots,u,i_1+i_2+\cdots+i_m=u$ , where $P_{i_1i_2
\ldots i_m}=P_{i_1i_2\ldots i_m}(x_1,x_2,\ldots x_m)$ are homogeneous
polynomials of degree $i_1$ relative to $x_1$, of degree $i_2$  relative  to
$x_2$ and so on, and finally, of degree $i_m$ relative to $x_m$, which
arise in computing of the power $(x_1+x_2+\cdots+x_m)^u$, and which are its
multihomogeneous components. On the other hand,
$$y(\omega_1x_1+\omega_2x_2+\cdots+\omega_m)^u=\sum\omega_1^{i_1}\omega_2
^{i_2}\ldots\omega_m^{i_m}yP_{i_1i_2\ldots i_m}(x_1,x_2,\ldots,x_m).$$
Hence $y(\omega_1x_1+\omega_2x_2+\cdots+\omega_mx_m)^u=0$.

Thus we have proved that the conditions  of the theorem hold in $R_\Omega=
(M_\Omega,L_\Omega)$. The first part of the proof, therefore, implies that
$L_\Omega$ is locally solvable.

{\sc LEMMA 5}. {\sl Let $R=(M,L)$ and $A$ is a finite-dimensional
subalgebra of $L$ such that $MA^n$ has finite dimension for some
integer $n$. Then a subspase of all $x\in L$ for which $xA^n=0$
has finite codimension}.

{\sc PROOF}. Let $z_{1},z_{2},\cdots,z_{k}$ is a basis of $A$.
Then any product of $n$ elements of the basis transfers $M$ into
finite-dimensional subspase of $M$. Thus the kernel of this
product has finite codimension. Since there is only a finite number
of different product of $n$ elements of the basis, the
intersection of all kernels of such products has a finite
codimension as well. The lemma is proved.

From Lemma 5 it follows that the associative algebra generated by
$L$ in $R=(M,L)$ is locally finite. Therefore we can apply the
results of [5] about the Jacobson radical of Lie algebra. Thus we
obtain the following

{\sc COROLLARY}. {\sl Let $L$ be a locally finite Lie algebra over a field
of characteristic 0. Suppose $L$ satisfies the condition  (3).  Then $L$ is
locally solvable if and only if $tr(ada)^2=0$  for every  $a\in L'$}.

{\sc PROOF}. The sufficiency of the condition is a consequence of Theorem 3,
since the kernel of the adjoint representation is the centre. Conversely,
assume $L$ is locally solvable. Then from [5] it follows that $ada, a\in
L'$, belongs to radical $J_{adL^*}(adL)$. Hence $ada$ is a nilpotent
linear transformation and $tr(ada)^2=0$.

\paragraph{3. A TRACE FORM.}

Let $R=(M,L)$ be a representation of a Lie algebra $L$ in a vector
space $M$ which satisfies condition (2). Then we can define a
trace form $f(a,b)=tra^Rb^R,a,b\in L$. The function $f(a,b)$ is
evidently a symmetric bilinear form on $M$ with values in the base
field $\Phi$. In particular, if $L$ satisfies (3), we obtain the
Killing form $K(a,b)=tr(ada)(adb)$.

If $f(a,b)$ is the trace form defined by the representation  $R=(M,L)$, then
$f([a,c],b)+f([a,[b, c])=tr([a, c]^Rb^R+a^R[b,c]^R)=tr([a^R,c^R]b^R+a^R[b^R,
c^R])=tr([a^Rb^R,c^R])=0$.

As noted in the Introduction, we can calculate the trace simultaneously for
every finite number of elements of $L$, since these elements can be
considered as linear transformations acting on a common quotient space of
M  of a finite dimension. Therefore, the last chain of equalities is correct.

A bilinear form $f(a,b)$ on $L$ that satisfies the condition
$f([a,c],b)+f(a,[b,c])=0$ is called an invariant form on $L$. Hence the
trace form is invariant. We note next that if $f(a,b)$ is any symmetric
invariant form on  $L$, then the radical $L^\perp$ of the form - that is,
the set of elements $z$ such that $f(a,z)=0$ for all $a\in L$ - is an
ideal. This is clear since $f(a,[z,b])=-f([a,b],z)=0$.

{\sc PROOF OF  THEOREM 4}. Let $f(a,b)$ be a trace form of $R=(M,L)$.
Then $L^\perp$ is an ideal of $L$ and $f(a,a)=tr(a^R)^2=0$ for every
$a\in L^\perp$. Hence $L^\perp$ is locally solvable by Theorem 3. Since $L$
is semi-simple, $L^\perp=0$, and $f(a,b)$ is non-degenerate. Next
suppose that the Killing form is non-degenerate. If $R(L)'\ne0$ then by
the Theorem 7 from [5] $R(L)'\subset J(L)$, Jacobson radical of $L$. This
implies that for every $a\in R(L)'$, it is true that $ada\in J(adL^*)$. Hence for
every $b\in L$ we have $ada\cdot adb\in J(adL^*)$. Since by Levitzki
theorem $J(adL^*)$ is locally nilpotent ideal, then $ada\cdot adb$ is a
nilpotent linear transformation and therefore $tr(ada\cdot adb)=0$.  This
contradicts our assumption that the trace form $tr(ada\cdot adb)$ is
non-degenerate. Therefore, it is required of $R(L)$ that $R(L)'=0$.

Let $a\in R(L)$, $b\in L$ and $N$ is a subspase of $L$ such that
$L /N$ is finite-dimensional and $ada$ and $adb$ act in $N$ as nil
transformations. Let us denote $\overline L=L/N$ and $\overline
{R(L)}=R(L)+N/N$ and choose a basis for $\overline L$ such that the
first vectors form a basis for $\overline{R(L)}$. The matrices of
linear transformations induced by $ada$ and $adb$ in $\overline
L$, respectively, are of the forms
$$\left(\begin{array}{cc}0&0\\**&0\end{array}
\right)\qquad\mbox{and}\qquad\left(\begin{array}{cc}*&0\\**&*\end{array}
\right).$$ This implies that $tr(ada)(adb)=0$. Hence $R(L)\subset
L^\perp$, and the Killing form is degenerate.

Using Killing form the following characterization of the locally
solvable radical in the characteristic 0 case (cf.\cite [p. 73]{6}) can
be obtained.

{\sc THEOREM 5}. {\sl If a locally finite Lie algebra $L$ over a field of
characteristic 0 satisfies the condition (3),  then the locally solvable
radical $R(L)$ of $L$ is the orthogonal complement $L'^\perp$ of $L'$
relative to the Killing  form $K(a,b)$}.

{\sc PROOF}. The algebra  $B=L'^\perp$ is an ideal. Further, if $b\in B'$, then
$tr(adb)^2=K(b,b)= 0$. The kernel of the representation  $a\to ada$, $a\in
B$, is abelian. Hence $B$ is locally solvable, by Corollary of theorem 2.
Hence $B\subset R(L)$. Next, let $x\in R(L)$, $a,b\in L$. Then $K(x,[a, b])
=K([x,a],b)$. By [5] the element $[x, a]$ belongs to the Jacobson
radical $J(L)$ of $L$. Consequently,  $ad[x,a]$  belongs to $J(adL^*)$ and
$ad[x,a]\cdot adb$ is nilpotent for every $b$.  Hence $K([x,a],b)=trad[x,
a]\cdot adb=0$. Thus $K(x,[a,b])=0$ and $x\in L'^\perp$. Thus $R(L)\subset
L'^\perp$ and so $R(L)=L'^\perp$.

Let $\Omega$ be an extension of the base field $\Phi$ of $L$. Then the
Killing form $f_\Omega$ of $L_\Omega$ is obtained from the Killing form $f$
of $L$ by the extension. Therefore,  $f_\Omega$ is non-degenerate if and
only if $f$ is non-degenerate \cite{8}. Consequently, $L_\Omega$ is
semi-simple if and only if $L$ is semi-simple (see also \cite{9}).

\paragraph{4. STRUCTURE  OF  SEMI-SIMPLE  ALGEBRAS.}

We continue the consideration of locally finite Lie algebras $L$
which satisfy the condition (3).

{\sc LEMMA 6} (cf. \cite[p. 29]{6}). {\sl Let  A  be a finite-dimensional
subalgera of $L$ such that $L(adA)^n\subset A$. Then
$A^\omega=\bigcap\limits_{k=1}^\infty A^k$ is an ideal of $L$}.

{\sc PROOF}. If $A^\omega=0$ then the assertion of Lemma is trivial. Let
$A^\omega\ne0$. We have $[L,A^n]\subset L(adA)^n$. $A$ is
finite-dimensional. Therefore, $A^\omega=A^m$ for some integer $m$. Thus $A^
\omega=A^{n+m-1}$. Now we have $[L,A^\omega]=[L,A^{n+m-1}]\subset L(adA)^{n+
m-1}\subset A(adA)^{m-1}=A^m=A^\omega$, which completes the proof.

Let $f(a,b)$ be any symmetric invariant bilinear form on $L$ and
let $A$ be a subspace of $L$. Denote by $A^\perp$ the subspace of
$L$ that consists of all elements $b\in L$ such that $f(a,b)=0$
for all $a\in A$.

{\sc LEMMA 7}. {\sl If  A  is a finite-dimensional subspace of $L$ such
that $A\cap A^\perp=0$, then $L=A\oplus A^\perp$}.

{\sc PROOF}. Since $A\cap A^\perp=0$, $A$ is a non-degenerate
subspace. Let us show that $L=A+A^\perp$. Take any basis
$a_1,a_2,\ldots,a_m$ in $A$ and let $c$ be any element of $L$.
Find a decomposition $c=a+a^\perp$, where $a \in A$ and
$a^\perp\in A^\perp$. We will look for $a$ in the form $a=x_1a_1+
x_2a_2+\cdots+x_ma_m$. Then $c$ will look like:
$c=x_1a_1+x_2a_2+\cdots+x_ma_m+a^\perp$. From $f(a_i,a^\perp)=0$
it follows that
$f(a_i,c)=\sum\limits_{k=1}^mx_kf(a_i,a_k),i=1,2,\ldots,m$. This
system of equations has exactly one solution, since its
determinant is the Gram's determinant of the system
$a_1,a_2,\ldots,a_m$. Since $A$ is a non-degenerate subspace, this
determinant is non-zero. The vector $a=
x_1a_1+x_2a_2+\cdots+x_ma_m$, where $x_k$ were just found, satisfies
the conditions $f(a_i,c-a)=0$. Indeed,
$f(a_i,c-a)=f(a_i,c)-\sum\limits_{k=1}^m x_mf(a_i,a_k)=0$. From
the equalities  $f(a_i,c-a)=0$, it follows that $c-a\in A^\perp$.
To complete the proof it remains to put $a^\perp=c-a$.

{\sc PROOF OF THEOREM 6}. Let $A$ be a finite-dimensional
subalgebra for which $L(adA)^n \subset A$. By Lemma 6, subalgebra
$A^\omega$ is an ideal of $L$. By Theorem 4,  the Killing form
$K(a, b)$ is non-degenerate. Since  $K(a,b)$ is invariant,
$(A^\omega)^\perp$ is an ideal of $L$. Indeed, for every $a\in
A^\omega$, every $b\in (A^\omega)^\perp$ and every $c\in L$ we
have  $K(a, [b,c])=-K([a,c],b)=0$. Let us prove that
$A^\omega\cap(A^\omega)^\perp=0$. If, on the contrary, $b_1,b_2\in
A^\omega\cap(A^\omega)^\perp$ and $a$ is any element of $L$, then
$K([b_1,b_2],a)=-K(b_1,[a,b_2])=0$. Since $K(a, b)$ is
non-degenerate, $[b_1,b_2]=0$. Hence
$A^\omega\cap(A^\omega)^\perp$ is an abelian ideal of $L$. Since
$L$ is semi-simple, $A^\omega\cap(A^\omega)^\perp=0$. By Lemma 7
this implies that $L=A^\omega\oplus (A^\omega)^\perp$. Since
$A^\omega$ and $(A^\omega)^\perp$ are direct summands, every ideal
of $A^\omega$ or $(A^\omega)^\perp$ is an ideal of $L$. Therefore
$A^\omega$ and $(A^\omega)^\perp$ are semi-simple subalgebras.
Moreover, since $A^\omega$ is finite-dimensional, $A^\omega$ is a
direct sum of ideals which are simple. Let us denote by $\Pi$ the
set of all finite-dimensional simple ideals of $L$. This set is
not empty; otherwise we have that $A^\omega=0$ for every
$x_{1},x_{2},\cdots,x_{k}\in L$. The previous implies that $L$ is
a locally nilpotent algebra. But by condition of the Theorem, $L$
is semi-simple. Since $M^\omega=M$ for every finite-dimensional
simple ideal $M$, $L=M\oplus M^\perp$. Now  denote by $N$ the
intersection of all $M^\perp$ where $M\in\Pi$. Let us prove that
$N=0$. Suppose not, and let $x_{1},x_{2},\cdots,x_{k}$ be a
non-zero element of $N$. Then by condition of the Theorem, there exists
a finite-dimensional subalgebra $A$, such that
$x_{1},x_{2},\cdots,x_{k}\in A$ and $L(adA)^n\subset A$. Let
$B=A\cap N$. Then $L(adB)^n\subset A$ and $L(adB)^n\subset N$,
since $N$ is an ideal of $L$. Consequently, $L(adB)^n\subset B$
and $B^\omega$ is ideal of $L$. If $B^\omega\ne0$, then $B^\omega$
is a semi-simple ideal contained in $N$. Next, $B^\omega$ is a
direct sum of simple ideals contained in $N$ and, consequently,
these ideals do not belong to $\Pi$. But this contradicts the
definition of $\Pi$. Hence $B^\omega=0$ for every finite subset
$x_{1},x_{2},\cdots,x_{k}\in N$ and $N$ is a locally nilpotent
ideal of $L$. But $L$ is a semi-simple Lie algebra. Hence $N=0$.
Then by Remac's theorem $L$ is a subdirect sum of simple
finite-dimensional ideals $M\approx L/M^\perp$. The proof of the
theorem is complete.

It is well known that every derivation of a semi-simple finite-dimensional
Lie algebra is inner. What about the infinite-dimensional case? Let $L =
\bigoplus\limits_\alpha L_\alpha$, where $L_\alpha$ are finite-dimensional
simple ideals of L. In this case Stewart~\cite{10} proved the theorem
that we give here for completeness.

{\sc THEOREM 7}. {\sl Let {L} be a semi-simple Lie algebra and let
$L$ is a direct sum, $L=\bigoplus \limits_\alpha L_\alpha$, where
$L_\alpha$ are finite-dimensional simple ideals of $L$. Then

1) every derivation of $L$ is an element of the complete direct sum $\tilde
L$ of $L_\alpha$ and the algebra of all derivations $D(L)$ of $L$ is
isomorphic to $\tilde L$.

2) a derivation $D$ is inner if and only if $L_\alpha D=0$ for every
$L_\alpha$ with the exception of a finite number of it

3) $L$ has outer derivations

4) every derivation of $L$ is locally finite.}

It follows from Theorem 7 that any semi-simple Lie algebra
satisfying condition (3) is a subalgebra of the algebra of all
derivations of such an algebra which is a direct sum of its simple
finite dimensional ideals.

\end{document}